\input amstex
\documentstyle{amsppt} 
\NoRunningHeads
\magnification=1150
\pageheight{7.25in}
\vcorrection{.25in}
\hcorrection{.25in}
\topmatter

\title Rigidity of area minimizing tori in 3-manifolds of nonnegative scalar curvature
\endtitle
\author Mingliang Cai and Gregory J. Galloway \endauthor
\affil Department of Mathematics and Computer Science \\ University of
Miami \endaffil
\address Coral Gables, FL 33124 \endaddress
\dedicatory  \enddedicatory
\thanks This work has been partially supported by NSF grant DMS-9803566.
\endthanks
\endtopmatter
\NoBlackBoxes
\document

\redefine\S{\Sigma}

\redefine \d{\partial}



\define\pd#1{\frac{\partial #1}{\partial u}}

\define\<{\langle}
\define\>{\rangle}

\define\g{\gamma}
\define\8{\infty}
\define\a{\alpha}

 The following conjecture arises from remarks in Fischer-Colbrie-Schoen
(\cite{FCS}, Remark 4, p. 207):
If $(M,g)$ is a complete Riemannian
$3$-manifold with nonnegative scalar curvature and if $\S$
is a two-sided torus in $M$  which is suitably of least area then $M$ is
flat.  Such a result,
as Fischer-Colbrie and Schoen commented, would be an interesting analogue
of the Cheeger-Gromoll
splitting theorem.  Here we present a proof of this conjecture assuming
$\S$ is of least area in its
isotopy class. The proof is a consequence of the following local
result, which is the main result of the paper.

\proclaim{Theorem 1} Let $(M,g)$ be a $C^{\8}$  $3$-manifold with
nonnegative scalar curvature,
$S \ge 0$. If $\S$ is a two-sided torus in $M$
which is locally of least area (see Section 2), then $M$ is  flat in a
neighborhood of $\S$.
\endproclaim

It follows that $\S$ is flat and totally geodesic and that locally $M$ splits
along
$\S$. A partial infinitesimal version of Theorem 1 was observed in
\cite{FCS}, namely, if $\S$ is a stable minimal two-sided torus in $M$
with nonnegative scalar curvature then $\S$ must be flat and totally
geodesic, and the
scalar curvature and normal Ricci curvature of $M$ vanishes along $\S$.  In
\cite{CG} the
authors proved Theorem~1 under the assumption that $M$ is analytic.
The result in the analytic case follows as
an immediate consequence of a more general result
which holds for $C^{\8}$ manifolds, see Theorem~B in \cite{CG}.
Here we will make use of Theorem B to present a proof of Theorem 1.

We note that, under the assumptions of Theorem 1, $M$ need not be globally
flat.
Consider, for example, $S^1\times S^2$, where $S^2$ is a sphere which is
flattened near the equator $E$.  Then $S^1\times E$ is a torus which is locally
of least area in $S^1\times S^2$.

The idea of the proof of Theorem 1 is as follows.  It is first shown that
$\S$ cannot be
locally {\it strictly\/} of least area.  If it were, then under a
sufficiently small
perturbation of the metric to a metric of (strictly) positive scalar
curvature, $\S$ would be
perturbed to a torus which is still locally of least area.  But this would
contradict
the fact that a compact two-sided stable minimal surface in a $3$-manifold
with strictly positive
scalar curvature must be a sphere, cf. Theorem 5.1 in
\cite{SY1}.  It is then shown that on each side of $\S$ there is a torus
which is locally
of least area.  By cutting out the region bounded by these two tori and pasting
it appropriately to a second copy, one obtains, using Theorem B, a smooth
$3$-torus with
nonnegative scalar curvature.  By Schoen and Yau \cite{SY1}, this $3$-torus
must be flat, and
Theorem~1 follows.  We now proceed to a detailed proof of Theorem 1.

In all that follows we work in the
$C^\8$ category.  For simplicity, all surfaces are assumed to be embedded.
However, by pulling
back to the normal bundle of $\S$, it is clear that a version of Theorem 1
holds for immersed
surfaces, as well.
By definition, a compact two-sided  surface $\S$ in a $3$-manifold $M$ is
locally of least area
provided in some normal neighborhood $V$ of $\S$, $A(\S) \le A(\S')$ for
all $\S'$ isotopic
to $\S$ in $V$, where $A $ is the area functional.  If the inequality is
strict for $\S'\ne\S$,
we say that
$\S$ is locally strictly of least area.
Note that ``locally of least area" implies ``stable minimal".

Let $V$ be a normal neighborhood of a compact two-sided surface $\S$ in a
$3$-manifold
$M$.  Then, via the normal exponential map, $V=(-\ell,\ell)\times \S$, and
the metric $g=ds^2$
takes the form,
$$
ds^2 = dt^2 + \sum_{i,j=1}^{2}g_{ij}(t,x)\, dx^idx^j \, . \tag1
$$

The following is a restatement of part of Theorem B in \cite{CG}.

\proclaim{Lemma 1} Let $(M,g)$ be a $3$-manifold with nonnegative scalar
curvature,
$S \ge 0$. Suppose $\S$ is a two-sided torus in $M$
which is locally of least area.  Then with respect to geodesic normal
coordinates
along $\S$ (see equation 1),
$$\frac{\d^n g_{ij}}{\d t^n}(0,x)= 0 \, ,
$$
for all $n$ and all $x\in\S$.
\endproclaim

Lemma 1 is used below to ensure that after certain cut and paste operations
the resulting
metric is smooth.

\proclaim{Lemma 2}  Suppose $\S$ is a compact two-sided surface in a
$3$-manifold $(M,g)$
with nonnegative scalar curvature, $S(g) \ge 0$.  Then there exists a
neighborhood $U$ of
$\S$ and a sequence of metrics $\{g_n\}$ on U such that $g_n \to g$ in
$C^{\8}$ topology on $U$,
and each $g_n$ has strictly positive scalar curvature, $S(g_n) > 0$.
\endproclaim

\demo{Proof} Let $V=(-\ell,\ell)\times\S$ be a normal neighborhood of $\S$,
so that the metric
$g$ takes the form (1).  Consider the sequence of metrics $g_n =
e^{-2n^{-1}t^2}g$.
A straight forward computation gives,
$$
S(g_n) =  e^{2n^{-1}t^2}(S(g) + 8n^{-1}(1+tH_t-n^{-1}t^2)),
$$
where $H_t$ is the mean curvature (in the metric $g$) of $\S_t=
\{t\}\times\S$.  It is clear
that by taking $\ell$
sufficiently small and $n$ sufficiently large we have $S(g_n) > 0$.
\enddemo

In the next lemma we show that if $\S\subset M $ is locally {\it
strictly\/} of least area then
by perturbing the metric of $M$ slightly,  $\S$ gets perturbed to a surface
which is still
locally of least area.

\proclaim{Lemma 3}  Suppose $\S$ is a compact two-sided surface in
$(M^3,g)$ which is locally
strictly of least area.  Let $\{g_n\}$ be any sequence of metrics such that
$g_n \to g$ in
$C^{\8}$ topology.  Then for any neighborhood $U$ of $\S$ and any positive
integer $N$
there exists, for some $n\ge N$, a surface $\S_n\subset U$ isotopic to $\S$
in $U$
which is locally of least area in $(M,g_n)$.
\endproclaim

\demo{Proof} The proof makes use of basic existence and convergence results
for least
area surfaces.  Let $V=[-\ell,\ell] \times \S$ be
a compact normal neighborhood of $\S$ contained in $U$, and restrict attention
to the compact Riemannian manifold-with-boundary $(V,g)$.
Since $\S$ is locally strictly of
least area, we can choose $\ell$ sufficiently small so that,
$$
A_g(\S) < A_g(\S') \quad\text{ for all } \S' \in \Cal I (\S), \,\S'\ne \S,
$$
where $\Cal I (\S)$ is the isotopy class of $\S$ in $V$, and $A_g$ is the
area functional
in the $g$ metric.

Set $V_0= [-\frac{\ell}2,\frac{\ell}2]\times\S$.
Let $f=f(t)$ be a smooth nonnegative function on $[-\ell,\ell]$ such that
$f=0$ on
$[-\frac{\ell}2,\frac{\ell}2]$.  By making the derivatives $f'(\pm \ell)$
sufficiently large in absolute
value, with $f'(\ell)>0$ and  $f'(-\ell)<0$,  we
obtain a conformally related metric $\bar g = e^fg$ with the following
properties.

\roster
\item $\bar g|_{V_0} = g|_{V_0}$.
\item $(V,\bar g)$  has strictly mean convex boundary, i.e., $\d V$ has
positive mean curvature.
\item For all  $\S'\in \Cal I(\S)$ such that $\S'\ne \S$, $A_{\bar g}(\S) <
A_{\bar g}(\S')$.
\endroster

For each $n$, set $\bar g_n = e^fg_n$.  Then the metrics $\bar g_n$
satisfy: (1)  $\bar
g_n|_{V_0} = g_n|_{V_0}$, (2)~$\bar g_n \to \bar g$ in $C^{\8}$ topology
and (3) for $n$
sufficiently large, $(V,\bar g_n)$ has mean convex boundary.  For each such
$n$ let,
$$
\a_n  =  \inf_{\S'\in \Cal I(\S)} A_{\bar g_n}(\S')\, .
$$
Then by Theorem 5.1 and Section 6 in \cite{HS} (see also \cite{MSY}) there
exists for each~
$n$ a surface $\S_n \in \Cal I(\S)$  such that $A_{\bar g_n}(\S_n) = \a_n$.
In applying the
results from  \cite{HS} we have used the fact that $V$ is $P^2$-irreducible
(provided $\S\ne S^2, P^2$) and that
$V$ does not contain any compact one-sided surfaces.  (If $\S = S^2$ or
$P^2$, one may appeal
to Theorem 5.2 in \cite{HS} and use specific features of the topology of
$[-\ell,\ell]\times
\S$).

For each $n$, $\S_n$ is a compact stable minimal surface in $(V,\bar g_n)$,
and the sequence
$\{\a_n\}$ is bounded.  It then follows by well-known convergence arguments
that there is
a subsequence of surfaces, call it again $\{\S_n\}$, which converges
locally in $C^{\8}$
topology to a compact (embedded) minimal surface $\bar \S$ in $(V,\bar g)$;
see especially
Section 2.2 in
\cite{M}, which applies fairly directly to the situation considered here.
By the nature of the
convergence, $\{\S_n\}$ is eventually contained in
any tubular neighborhood of $\bar \S$, and for $n$
sufficiently large, $\S_n$ will be transverse to the normal geodesics of
$\bar\S$.  It follows
that $\S_n$ covers $\bar \S$ via projection along the normal geodesics.
Since $\bar \S$ is
necessarily two-sided (again, because $V$ does not contain any compact
one-sided surfaces),
it follows that the covering of $\bar \S$ by $\S_n$ must be one-sheeted,
i.e., projection along
the normal geodesics of $\bar\S$ provides a diffeomorphism of $\S_n$ onto
$\bar \S$; see e.g.,
\cite{S}.

Thus, $\bar\S$ is isotopic to $\S$ since each $\S_n$ is.
Furthermore, we have,
$$
A_{\bar g}(\bar \S) = \lim_{n\to\8} \a_n \le \lim_{n\to\8} A_{\bar g_n}(\S)
= A_{\bar g}(\S)\, .
$$
Since $\S$ is strictly of least area in its isotopy class in $(V,\bar g)$,
we conclude that $\bar
\S =
\S$. But by the above convergence, this means that for $n$ large enough,
$\S_n$ is contained in
$\operatorname{int}V_0$, in which $\bar g_n = g_n$.  It follows that, for
such $n$, $\S_n$
is locally of least area in $(V,g_n)$.  This concludes the proof of Lemma
3.
\enddemo

\demo{Proof of Theorem 1} Let $V = (-\ell,\ell)\times \S$ be a normal
neighborhood
 of $\S$ with metric $g$ as in equation (1).  Choose $\ell$ sufficiently
small so that
$\S$ is of least area in its isotopy class in V.  For technical reasons we
modify the metric
$g$ as follows.  Let $\hat g$ be the metric on $V$ of the form (1) but with
component functions
$\hat g_{ij}$, $1\le i,j \le 2$, defined by,
$$
\hat g_{ij}(t,x) = \cases g_{ij}(t,x),&\text{for $t\in[0,\ell)$}\\
g_{ij}(-t,x),&\text{for $t\in (-\ell,0]$}. \endcases
$$

$(V,\hat g)$ is a smooth Riemannian manifold (by Lemma 1) such that $S(\hat
g) \ge 0$
and reflection across
$\S$, $(t,x) \mapsto (-t,x)$, is an isometry.  Further, $\S$ is of least
area in its
isotopy class in $(V,\hat g)$.  By choosing $\ell$ even smaller if necessary,
we guarantee that Lemma~2 holds for the neighborhood $U=V$.

If $\S$ were {\it strictly\/} of least area in its isotopy class in
$(V,\hat g)$ then Lemmas 2
and 3 would imply that there is a two-sided stable minimal torus $\S'$ near
$\S$ with respect
to some metric of strictly positive scalar curvature on $V$. This would
contradict
Theorem 5.1 in \cite{SY1}.  Thus, there exists a surface $\bar \S \in \Cal
I(\S)$, $\bar \S
\ne\S$ such that $A_{\hat g}(\bar \S) = A_{\hat g}(\S)$.  Hence, $\bar \S$
is also of
least area in its isotopy class.

We claim that $\bar \S$ is contained in one of the components of
$V\setminus \S$.
If not, then $\bar \S$ and $\S$ must meet.  Since $\bar \S$ and $\S$ are
stable minimal tori
in $(V,\hat g)$ they must be totally geodesic (cf. \cite{FCS}).  Since they
are totally
geodesic and distinct, they must meet transversally.  Thus, the
intersection of $\bar\S$ and
$\S$ will consist of a finite number of circles.  By reflecting the portion
of $\bar \S$
in $(-\ell,0]\times \S$ across $\S$ to $[0,\ell)\times\S$ and smoothing out
the resulting ridge
along the circles of intersection, we obtain a surface isotopic to $\S$
with less area than
$\S$, which is a contradiction.  Thus, $\bar \S $ lies to one side of $\S$
and does not meet
$\S$.

These arguments imply that in the original Riemannian manifold $(V,g)$
there exist two
tori $\S^+$ and $\S^-$, one on each side of $\S$,  each isotopic to $\S$
and each locally
of least area.  Let $W$ be the region in $V$ bounded by  $\S^+$ and $\S^-$.
Standard properties of isotopies \cite{H} imply that
$W$ has topology $[-1,1]\times T^2$.
By taking two copies of $W$ and gluing them appropriately along their
boundaries, we obtain, by Lemma 1, a smooth Riemannian manifold with
nonnegative scalar
curvature which is diffeomorphic to a $3$-torus.  By Schoen and Yau \cite{SY1},
this $3$-torus, and hence $W$ must be flat.
\enddemo

By fairly standard continuation arguments, Theorem 1 can be globalized as
follows.

\proclaim{Theorem 2} Let $M$ be a  complete connected $3$-manifold of
nonnegative scalar
curvature whose boundary (possibly empty) is mean convex.  If  $M$ contains
a two-sided
torus $\S$ which is of least area in its isotopy class then $M$ is flat.
\endproclaim

\demo{Proof} By the maximum principle, either $\S$ is a boundary component
of $M$
or $\S$ is in the interior of $M$.  If $\S$ is a boundary component, let
$M_0=M$.
If $\S$ is in the interior and disconnects $M$, let $M_0= \overline U_0$,
where $U_0$ is one of the
components of $M\setminus \S$.  If $\S$ is in the interior and does not
disconnect $M$, let $M_0$ be the
manifold  with boundary obtained by
``separating" $M$ along $\S$.
In all cases, $\S$ is a boundary component of $M_0$. To prove Theorem 2 it
suffices to show that $M_0$ is flat.  Consider the normal exponential map
$\Phi: [0,\8)\times\S \to M_0$ along $\S$ defined by $\Phi(t,x)= \exp_xtN$,
where N is the
inward pointing unit normal along $\S$.  (Note $\Phi$ need not be defined
on all of
$[0,\8)\times\S$.)

By Theorem 1, $M_0$ is flat in a neighborhood of $\S$.  (It is easily seen that
Theorem~1 is still valid if $\S$ is a boundary component.)
Then, by standard arguments (which
require only nonnegative Ricci curvature), since  $\S$ is locally of
least area  there exists $a>0$ such that $\Phi: [0,a)\times\S \to
\Phi([0,a)\times\S)$ is an isometry.
(Here $[0,a)\times\S$ carries the standard product metric and hence is flat
since $\S$ is).  Let $\ell$
be the largest number (possibly $\8$) such that $\Phi: [0,\ell)\times\S \to
\Phi([0,\ell)\times\S)$
is an isometry.  Consider first the case $\ell =\8$. Using that the limit
of a sequence of normal geodesics
to $\S$ is a normal geodesic to $\S$, one easily verifies that
$\Phi([0,\8)\times\S)$
is both open and closed in $M_0$.  Hence, $\Phi([0,\8)\times\S)=M_0$ and
$M_0$ is flat.

Now consider the case $\ell < \8$.  Since $M_0$ is complete, each normal
geodesic to $\S$, $\g_x: t
\mapsto \Phi(t,x)$, $0\le t< \ell$, extends to $t=\ell$.  Suppose that
$\Phi: [0,\ell]\times\S \to
\Phi([0,\ell]\times\S)$ is an isometry.  Then $\S_{\ell}=
\Phi(\{\ell\}\times \S)$ is an embedded
totally geodesic torus in $M_0$ which is locally of least area.  By the
maximality of $\ell$, $\S_{\ell}$ must
meet some component
$\S'$ of $\d M_0$.  By the maximum principle for hypersurfaces, $\S_{\ell}$
and $\S'$ must agree,
$\S_{\ell}=\S'$.  One can now argue that $\Phi([0,\ell]\times\S)$ is both
open and closed in
$M_0$. Hence, $M_0=\Phi([0,\ell]\times\S)$ is flat.

Now suppose $\Phi: [0,\ell]\times\S \to \Phi([0,\ell]\times\S)$ is not an
isometry.  The only way this
can happen is if two normal geodesics to $\S$, $\g_{x_i}$, $i=1,2$, meet at
$t=\ell$,
$\g_{x_1}(\ell)=\g_{x_2}(\ell)$.   Since there can be no focal points to
$\S$ along $\g_{x_i}|_{[0,\ell]}$,
there exists a neighborhood $U_{i}$ of $x_i$ in $\S$ such that $\Phi:
[0,\ell]\times U_i \to
\Phi([0,\ell]\times U_i)$ is an isometry.  Hence, $\Phi(\{\ell\}\times
U_i)$ is an embedded totally
geodesic hypersurface in $M_0$ which (by the choice of $\ell$) is a
constant distance $\ell$ from $\S$.  It
follows that
$\Phi(\{\ell\}\times U_1)$ and $\Phi(\{\ell\}\times U_2)$ must agree near
the common
end point $\g_{x_1}(\ell)=\g_{x_2}(\ell)$.
By a straight forward continuation argument we conclude that each geodesic
segment  $\g_x$, $x\in\S$,
of length $2\,\ell$ meets $\S$ orthogonally at both end points.  It is now
easily argued
that $\Phi([0,\ell]\times\S)$ is both open and closed in $M_0$.  Hence,
$M_0=\Phi([0,\ell]\times\S)$ is flat,
and the proof of Theorem~2 is complete.
\enddemo

We make some concluding remarks.
\smallskip
\noindent
1. The results presented here were motivated in part by certain problems
concerning the topology of black holes in general relativity, cf., \cite{CG}
\cite{G1}, \cite{G2}.

\smallskip
\noindent
2.  The example mentioned after Theorem 1,  $M = S^1\times S^2$, with $S^2$
flattened near the equator $E$,
and $\S= S^1\times E$, shows that stability is not sufficient to imply
flatness.  Assume the $S^1$ factor and $E$ have the same radius.
Cutting $M$ along $\S$ we obtain two solid tori, the boundary of each of
which is a copy of $\S$.
Gluing the two solid tori back together along their toroidal boundaries
after a suitable twist
we obtain a manifold $M'$ diffeomorphic to $S^3$ with nonnegative scalar curvature  
 which contains a stable minimal torus.
Applying Theorem A in \cite{GL}, which is proved by a local construction,
we can add an
asymptotically flat end to $M'$ to obtain an asymptotically flat manifold
diffeomorphic to $\Bbb R^3$ with nonnegative scalar curvature which contains
a stable
minimal torus.  In the language of general relativity, we have obtained an
asymptotcally flat
time symmetric initial data set on $\Bbb R^3$ satisfying the constraint
equations
which contains a stable toroidal apparent horizon.  However, we know of no
such
vacuum (scalar flat) examples, and conjecture that there are none.

\smallskip
\noindent
3.  Using the higher dimensional
work of Schoen and Yau \cite{SY2} it appears that Theorem 1 can be extended
to higher dimensions
as follows: Let $M^n$ have nonnegative scalar curvature. If $\S$ is a
compact two-sided
hypersurface in $M^n$ which does not admit a metric of positive scalar
curvature and which is locally
of least area then a neighborhood of $\S$ splits.  We are grateful to a
referee for
suggesting an alterative proof of Lemma 3, valid in higher dimensions,
which makes
this generalization possible.  Further aspects of this  will be discussed
elsewhere.

\smallskip
\noindent
4. In \cite{FCS}, Fischer-Colbrie and Schoen proved that a complete stable
minimal surface
in an orientable $3$-manifold with nonnegative scalar curvature must be
conformal to the
complex plane or the cylinder $A$.  In the latter case it has been shown that $A$
is flat and totally geodesic, cf., \cite{FCS} and \cite{CM}.  The example $M=\Bbb R\times S^2$,
where $S^2$
is flattened near the equator, shows that $M$ need not be flat.  However,
in view of the
results cited and the results presented here, it seems reasonable to conjecture
that if the cylinder $A$ is actually area minimizing (in a suitable sense)
then $M$ is flat
(cf., Remark 4 in \cite{FCS}).

\smallskip
We would like to express our gratitude to a referee and an anonymous
reviewer for many valuable
comments and suggestions.  We would also like to thank Tobias Colding for
some helpful
discussions.

\enddocument